\documentclass[aps,twocolumn,showpacs]{revtex4}
\usepackage[dvips]{graphicx}

\begin{document}
\title{New power-law tailed distributions emerging in $\kappa$-statistics}
\date{\today}
\author{G. Kaniadakis}
\email{giorgio.kaniadakis@polito.it}
\affiliation{Department of Applied Science and Technology, Politecnico di Torino, \\
Corso Duca degli Abruzzi 24, 10129 Torino, Italy}
\begin{abstract} 
Over the last two decades, it has been argued that the Lorentz transformation mechanism, which imposes the generalization of Newton's classical mechanics into Einstein's special relativity, implies a generalization, or deformation, of the ordinary statistical mechanics. The exponential function, which defines the Boltzmann's factor, emerges properly deformed within this formalism. Starting from this, so-called $\kappa$-deformed exponential function, we introduce new classes of statistical distributions emerging as the $\kappa$-deformed version of already known distribution as the Generalized Gamma, Weibull, Logistic which can be adopted in the analysis of statistical data that exhibit power-law tails.
\end{abstract}

\maketitle


\section{Introduction}

There is a vast phenomenology related to experimental statistical distributions with a bulk described by well known theoretical exponential models, in the following called classical models, like the Generalized Gamma distribution, the Weibull distribution, the Logistic distribution etc, while it seems that their tails are gradually substituted by fat, not exponential tails following Pareto's law $p(x)= A\,x^{-b}$, \cite{PowerLawKaniadakis} (and references therein). This anomaly regards not only physical systems but also other natural or artificial systems.

Let us focus on a specific classical model as for instance the Generalized Gamma distribution though discussion can be extended also to Weibull and Logistic models. It is important to emphasize that the Gamma model and the Pareto model refer to two distinct spectral regions, namely the lower and the upper region of the spectrum respectively. This dichotomy gives rise to the question about whether the two distributions are actually two approximations, in the low and in the high spectral regions respectively, of a unique theoretical distribution, which hold over in the whole spectrum. A further question that arises spontaneously is whether there exists a simple and transparent underlying mechanism, based on first principle, that generates this unique distribution.

With reference to physical systems, the metaphor of relativistic physics is of great help in the search for the mechanism that generates  a statistical distribution holding at both low and high energies.
It should be recalled that the laws of relativistic physics and the involved physical quantities emerge as generalizations of the corresponding ones of classical physics. For instance, the dimensionless relativistic kinetic energy is given by ${\cal E}_{\kappa}=\left (  \sqrt{1+\kappa^2 q^2}-1 \right )/{\kappa^2}$, where $q$ is the dimensionless momentum and $\kappa$ is the reciprocal of the light speed, written in dimensionless form. This relativistic expression of kinetic energy can be viewed as a generalization or deformation of classical kinetic energy through the deformation parameter $0<\kappa<1$. In the classical limit, corresponding to $q \rightarrow 0$,  or alternatively to $\kappa \rightarrow 0$,  the classical expression of kinetic energy ${\cal E}_{0}\approx q^2/2$ is obtained while in the ultra relativistic limit, corresponding to $q \rightarrow +\infty$,  ${\cal E}_{\kappa}\approx q/\kappa$ is obtained. When one moves from classical to relativistic physics, not only the expressions of the various physical quantities (momentum, velocity, total energy, force etc) but also all the mathematical tools of the theory (relativistic additivity law of velocities, momenta additivity law, Lorentz invariant integral, relativistic derivative etc) emerge properly generalized or deformed. This $\kappa$-deformed mathematical formalism of special relativity necessarily leads to a proper deformation of all the mathematical functions and in particular of the exponential function which expression follows directly from the energy-momentum Lorentz transformations \cite{PRE2002} and assumes the form
\begin{eqnarray}
\exp_{\kappa}(x)= \left (\sqrt{1+\kappa ^2 x^2}+\kappa x \right ) ^{1/\kappa} \ .
\label{KGGDI3}
\end{eqnarray}

This $\kappa$-deformed exponential, or simply $\kappa-$exponential is defined over the whole real axis, and its most interesting feature is undoubtedly given by its asymptotic behaviour. In the low energy limit, obtained when $q \rightarrow 0$ or equivalently when $\kappa \rightarrow 0$, the $\kappa$-exponential reduces to the Euler ordinary exponential i.e. $\exp_{\kappa}(x)\,\,
{\atop\stackrel{\textstyle\approx}{\scriptstyle x\rightarrow
0}} \,\,\, \exp\,(x)$.  On the other hand, in the high energy limit, the $\kappa$-exponential exhibits power-law tails described by Pareto's law i.e. $\exp_{\kappa}(x)\,\,
{\atop\stackrel{\textstyle\approx}{\scriptstyle x\rightarrow
\pm \infty}} \,\,\,|2 \kappa \, x |^{\pm 1/\kappa}$.

The $\kappa$-exponential represents a very useful and
powerful tool to formulate a generalized statistical theory capable to treat systems described by distribution functions exhibiting power-law
tails \cite{PRE2002,PRE2017,SciRep2020}. Generalized statistical mechanics, based on $\kappa$-exponential, preserves the main features of the Boltzmann-Gibbs statistical mechanics based on the ordinary exponential through the Boltzmann factor.

Over the last two decades, the $\kappa$-statistical theory has attracted the interest of many researchers, who have studied its foundations \cite{Silva06A,Naudts1,Topsoe,Scarfone2013,SouzaPLA2014,daSilva2020A}, and the underlying thermodynamics \cite{Wada1,ScarfoneWada,Bento3lawThermod,WadaMatsuzoeScarfone2015,Mehri2020}, and at the same time, have considered specific applications of the theory to various fields of science.  A nonexhaustive list of applications includes, among others, those in quantum statistics \cite{Santos2011a,Planck,Lourek2}, in quantum theory \cite{Ourabah,Costa2020,Andrade2020}, in plasma physics \cite{Lourek,Gougam,Chen,Landau2017,Qualitative2017,Lourek2019,Kalid2020A}, in nuclear fission \cite{NuclearEnergy2017,Guedesa2021}, in particle physics \cite{Shen2020},  in astrophysics \cite{Carvalho,Carvalho2,Carvalho2010,Cure,Soares2019}, in cosmology \cite{AbreuEPL,AbreuIJMPA,ChinesePL,AbreuEPL2018A,Immirzi,Yang2020,Moradpour2020}, in geomechanics \cite{Oreste2017,Oreste2019}, in genomics \cite{SouzaEPL2014,Costa2019}, in complex networks \cite{Macedo,Stella}, in waveform inversion algorithms \cite{daSilva2020}, in image processing \cite{Lei2020}, in machine learning \cite{Passos2020}, in seismology \cite{seismos},
 in economy \cite{Clementi2010,Modanese,Vallejos2019} and in finance \cite{Trivellato2012,Trivellato2013,Tapiero,Trivellato2017}.

Main goal of the present effort is the proposal of five classes of statistical distribution with support $x>0$ (but easily generalized to the case $x \in R$) presenting power-law tails, starting from the above mentioned $\kappa$-exponential function. For all these distributions their momenta of any order are obtained in closed forms so that various properties like the mean, the variance, the coefficient of variation, the skewness and the kurtosis can be evaluated easily.

Table 1 shows  the correspondence between the proposed distributions and their already known ordinary counterparts (if applicable), obtained in the classical $\kappa \rightarrow 0$ limit.

\begin{table}[h]
	\caption{Correspondence between the here proposed distributions and their already known ordinary counterpart\\}
    \centering
	\begin{tabular}{l c c} \hline
		\textbf{Name} & \textbf{Corresponding known distribution} & \textbf{} \\ \hline
		Type I	& Generalized Gamma distribution& \\
		Type II &	Weibull distribution &	 \\
		Type III &	Generalized Logistic distribution &	 \\
		Type IV	& not defined &  \\
		Type V &Exponential distribution&  \\
        \hline
	\end{tabular}
\end{table}

The present approach can be easily extended in order to include the $\kappa$-deformed versions of all the statistical distributions, already known in the literature, involving the ordinary exponential function.

\section{Distribution Functions of Type I}

Let us consider the cdf $P_{\kappa}(x)$ defined for $x\geq 0$, through
\begin{eqnarray}
P_{\kappa}(x)=\frac{\int_{\,0}^{x}h(t)
\,\exp_{\kappa}[-f(t)]\,dt}{\int_{\,0}^{\infty}h(t)
\,\exp_{\kappa}[-f(t)]\,dt} \ \ , \label{MSTI0A1}
\end{eqnarray}
where the functions $g(x)$ and $f(x)$ are given by
\begin{eqnarray}
&&f(x)=\beta x^{\alpha} \ \ , \label{MSTI0A4} \\
&&h(x)=x^{\alpha \nu-1} \ \ , \label{MSTI0A3}
\end{eqnarray}
with $\alpha > 0$, $\beta>0$ and $0 < \nu < 1/ \kappa$. $P_{\kappa}(x)$ is a non-decreasing function of $x$ taking values in the interval $[0, 1]$.

After performing the integration in the denominator appearing in Eq. (\ref{MSTI0A1}), the cdf assumes the form
\begin{eqnarray}
P_{\kappa}(x)=N_{\kappa}
\,\, \int_{\,
0}^{ x}\!\!z^{\alpha \nu-1} \,\exp_{\kappa}(-\beta z^{\alpha})\,dz \ \ , \label{KGGDII5}
\end{eqnarray}
where the normalization constant \cite{SupplementaryMaterial} is given by
\begin{eqnarray}
N_{\kappa}=(1+\kappa \,\nu)\,(2\,\kappa)^{\nu}
 \,\,
\frac{\Gamma\!\left(\frac{1}{2\kappa}+\frac{\nu}{2}
\right)}{\Gamma\!\left(\frac{1}{2\kappa}-\frac{\nu}{2}
\right)}\,\, \frac{\alpha\,\beta^{\,\nu}}{\Gamma\!\left(\nu\right)} \ \ .\label{KGGDII6}
\end{eqnarray}

The pdf defined by means of $p_{\kappa}(x)=\frac{d \, P_{\kappa}(x)}{dx}$, i.e.
\begin{eqnarray}
p_{\,\kappa}\,(x)= N_{\kappa}
\,\, x^{\alpha\nu-1} \,\exp_{\kappa}(-\beta x^{\alpha}) \ \ , \label{KGGDII7}
\end{eqnarray}
represents the $\kappa$-deformed version of the Generalized Gamma  pdf and in the $\kappa \rightarrow 0$ limit reduces to the ordinary Generalized Gamma pdf $p(x)=\frac{|\alpha|\,\beta^{\,\nu}}{\Gamma(\nu)}\,\,x^{\alpha\nu-1} \exp(-\beta x^{ \alpha})$.

Asymptotically the pdf $p_{\,\kappa}(x)$, behaves according to
\begin{eqnarray}
&&p_{\,\kappa}(x)\,\,
{\atop\stackrel{\textstyle\sim}{\scriptstyle x\rightarrow 0^+}}
\,\,\,N_{\kappa}\,\,\, x^{\alpha\nu-1} \ \ . \label{KGGDII7a} \\
&&p_{\,\kappa}(x)\,\,
{\atop\stackrel{\textstyle\sim}{\scriptstyle x\rightarrow +\infty}}
\,\,\,(2\,\kappa\,\beta)^{\!-1/\kappa}\,\,N_{\kappa}\,\, x^{\alpha\nu-1-\alpha/\kappa} \ \ . \label{KGGDII7b}
\end{eqnarray}

For $0<\nu < \frac{1}{\alpha} < \frac{1}{\kappa} $, $p_{\kappa}(x)$ decreases monotonically from the $p_{\kappa}(0^+)=+\infty$ value to the $p_{\kappa}(+\infty)=0^+$ value, and therefore behaves qualitatively like the Pareto distribution. The main difference, with respect the Pareto distribution, is that $p_{\kappa}(x)$ has a unitary norm over the entire non-negative real axis i.e. $\int_{0}^{+\infty}p(x)dx=1$.
For $0<\nu = \frac{1}{\alpha} < \frac{1}{\kappa} $, $p_{\kappa}(x)$ decreases monotonically  from the $p_{\kappa}(0^+)=N_{\kappa}$ value to the $p_{\kappa}(+\infty)=0^+$ value. Finally, for $0<\frac{1}{\alpha} <  \nu  < \frac{1}{\kappa} $, $p_{\kappa}(x)$ starts by the zero value and in
\begin{eqnarray}
x_{mode}=\beta^{-1/\alpha}\, \left (\nu -\frac{1}{\alpha} \right )^{\frac{1}{\alpha}}
\left [1-\kappa^2\left(\nu-\frac{1}{\alpha}\right)^{\!2}\, \right ]^{-\frac{1}{2\alpha}}  , \label{KGGDII10}
\end{eqnarray}
presents its maximum value and after decreases monotonically to the $p_{\kappa}(+\infty)=0^+$ value.

The moment of order $m$ of the pdf $p_{\,\kappa}(x)$, defined by means of $<x^m>\,=\int_0^{+\infty}\!\!x^m\, p_{\,\kappa}(x)\, dx$, has a finite value if $0<\nu+m/\alpha
<1/\kappa$, and can be expressed in  closed form, in terms of the Euler $\Gamma(x)$ function \cite{SupplementaryMaterial}, through
\begin{eqnarray}
<x^m>\,=\!\!\!\!\!\!\!\!\!\!&&\beta^{-m/\alpha}\,\frac{(1+\kappa\nu)\,(2\kappa)^{-m/\alpha}}
{1+\kappa\left(\nu+\frac{m}{\alpha}\right)}\,\frac{\Gamma\left(\nu+\frac{m}{\alpha}\right)}{\Gamma\left(\nu\right)} \nonumber \\
\!\!\!\!\!\!\!\!\!\!&&\times\,\frac{\Gamma\left(\frac{1}{2\kappa}+\frac{\nu}{2}\right)}
{\Gamma\left(\frac{1}{2\kappa}-\frac{\nu}{2}\right)}\,\,
\frac{\Gamma\left(\frac{1}{2\kappa}-\frac{\nu}{2}-\frac{m}{2\alpha}\right)}
{\Gamma\left(\frac{1}{2\kappa}+\frac{\nu}{2}+\frac{m}{2\alpha}\right)}
\ \ .  \label{KGGDII11M}
\end{eqnarray}

In Table 2 are reported the correspondences of some special cases of pdf defined in Eq. (\ref{KGGDII7}) with their classical counterpart.

\begin{table}[h]
	\caption{Special cases of the distribution of type I and their correspondences with already known distributions which recover in the $\kappa \rightarrow 0$ limit \\}
    \centering
	\begin{tabular}{l c c} \hline
		\textbf{Name} & \textbf{Corresponding known distribution} & \textbf{} \\ \hline
		Type $I_1$ &Exponential ($\alpha=1$, $\nu=1$) & \\
		Type $I_2$ &Erlang  ($\alpha=1$, $\nu=$integer) &	 \\
		Type $I_3$ &Gamma  ($\alpha=1$) &	 \\
		Type $I_4$ &Chi-Squared ($\alpha=1$, $\nu=$ half integer) &  \\
        Type $I_5$ & Nakagami ($\alpha=2$) $\nu>0$) & \\
		Type $I_6$ &Rayleigh ($\alpha=2$, $\nu=1$) &	 \\
		Type $I_7$ &Chi ($\alpha=2$, $\nu=$ half integer)&	 \\
        Type $I_8$ &Maxwell ($\alpha=2$, $\nu=3/2$)&  \\
        Type $I_9$ &Half-Normal ($\alpha=2$, $\nu=1/2$)& \\
		Type $I_{10}$ &Weibull ($\alpha>0$, $\nu=1$)&	 \\
		Type $I_{11}$ &Stretched Exponential ($\alpha>0$, $\nu=1/\alpha$)&	 \\
        \hline
	\end{tabular}
\end{table}

\subsection{The distribution of Type I$_1$}

The special case of the distribution of Type I corresponding ${\alpha}=\nu=1$ defines the distribution of Type I$_1$. In this case the integral in the expression of cdf can be computed easily so that the pdf $p_{\,\kappa}(x)$ and the cdf $P_{\,\kappa}(x)$ assume the following simple form
\begin{eqnarray}
\!\!\!\!\!\!\!\!\!\!\!\!&&p_{\,\kappa}(x)=(1-\kappa^2)\,\beta\,\exp_{\kappa}(-\beta \,x) \ \ ,
\\
\!\!\!\!\!\!\!\!\!\!\!\!&&P_{\,\kappa}(x)= 1 -
\left(\sqrt{1+\kappa^2\beta^2x^2}+\kappa^2\beta x \right)
\exp_{\kappa}(-\beta x)\ \ . \label{MSTI16} \ \ \ \ \ \ \ \ \
\end{eqnarray}
and represent the $\kappa$-deformed version of the Exponential distribution which recovers in the $\kappa \rightarrow 0$ limit.

The moment of order $m$ of $p_{\,\kappa}(x)$ has a finite value if $0<m +1 <1/\kappa$ and is given by
\begin{eqnarray}
<x^m>\,=\displaystyle{\frac{1-\kappa^2}{
\prod_{n=0}^{m+1}[1-(2n-m-1)\,\kappa]}\,\,\frac{m!}{\beta^{m}}} \ \ . \label{MSTI12}
\end{eqnarray}
The expectation value and variance of $p_{\,\kappa}(x)$ are
given by
\begin{eqnarray}
\!\!\!\!\!&&<x>\,= \frac{1}{\beta}\,\frac{1-\kappa^2}{1-4\kappa^2} \
\ , \label{MSTI13} \\ \!\!\!\!\!&&\sigma^2=
\frac{1}{{\beta}^{2}}\,\frac{2(1-4\kappa^2)^2 - (1-\kappa^2)^2
(1-9\kappa^2)}{(1-4\kappa^2)^2(1-9\kappa^2)}   \ \ . \ \ \ \ \ \
\label{MSTI14}
\end{eqnarray}

\subsection{The distribution of Type I$_2$}

The most general subclass of distributions of Type I, whose cdf can be obtaining in closed forms, in following called distributions of Type I$_2$, corresponds to the choice $\alpha=1$ and $\nu=n=$positive integer. These distributions emerges to be the $\kappa$-deformed version of the ordinary Erlang distributions. The first element of the family of distributions I$_2$ is just the distribution I$_1$. An important property of distributions I$_2$ is that the corresponding  survival functions $S_{\kappa}(x)=1-P_{\kappa}(x)$ can be expressed in closed form and therefore the rate equation $d\,S_{\kappa}(x)/dx=-h_{\kappa}(x)\, S_{\kappa}(x)$  permits to obtain in closed form also the hazard functions $h_{\kappa}(x)$.

Hereafter we propose a procedure to compute $P_{\kappa}(x)$ of Type I$_2$ starting from the related $p_{\kappa}(x)$ after posing, for simplicity of the exposition, $\beta=1$ i.e.
\begin{equation}
p_{\kappa}(x)=N_{\kappa}\,\,x^{n-1} \exp_{\kappa}(-x)  \ \ , \label{KGGDIII1}
\end{equation}
$n$ being a positive integer and
\begin{equation}
N_{\kappa}= \frac{1}{(n-1)!} \prod_{m=0}^n [1+(2m-n)\kappa]  \, \,
  \ \ . \label{KGGDIII2}
\end{equation}
The $P_{\kappa}(x)$ assumes the form
\begin{equation}
P_{\kappa}(x)=N_{\kappa}\int_0^{x}t^{n-1} \exp_{\kappa}(-x)  \ \ , \label{KGGDIII3}
\end{equation}
and in order to perform the integral that appears in Eq. (\ref{KGGDIII3}), we introduce the ansatz
\begin{equation}
P_{\kappa}(x)=1 - \left [ \,R_{\kappa}(x) + Q_{\kappa}(x) \sqrt{1+\kappa^2 x^2 }\,\,\right ] \exp_{\kappa}(-x)  \ \ , \label{KGGDIII4}
\end{equation}
where $R_{\kappa}=R_{\kappa}(x)$ and $Q_{\kappa}=Q_{\kappa}(x)$ are two unknown functions that have to be determined. After substitution of this ansatz in the left hand side of Eq. (\ref{KGGDIII3}), and after derivation of both sides of the equation, it obtains
\begin{equation}
A + B \, \sqrt{1+\kappa^2 x^2 }=0 \ \ . \label{KGGDIII5}
\end{equation}
with
\begin{eqnarray}
&&A= (1+\kappa^2 x^2) \frac{d Q_{\kappa}}{dx} +\kappa^2 x\, Q_{\kappa}- R_{\kappa} \ \ , \label{KGGDIII6A}
\\
&&B=\frac{d R_{\kappa}}{dx} - Q_{\kappa} + N_{\kappa} x^{n-1} \ \ . \label{KGGDIII6B}
\end{eqnarray}

Eq. (\ref{KGGDIII5}) holds for any value of the variable $x$, and this is guaranteed if and only if $A=0$ and $B=0$. These two conditions represent a system of two coupled first-order differential equations for the unknown functions $R_{\kappa}$ and $Q_{\kappa}$. After uncoupling this system it obtains that the function $R_{\kappa}$ obeys the following  differential equation
\begin{eqnarray}
&&(1+\kappa^2 x^2)\,\frac{d^2 R_{\kappa}(x)}{d\,x^2} + \kappa^2 x \, \frac{d R_{\kappa}(x)}{d\,x} - R_{\kappa}(x) \nonumber \\ && + \kappa^2\,n\, N_{\kappa}\, x ^{n} + (n-1)\, N_{\kappa}\, x ^{n-2}=0 \ \ , \label{KGGDIII7}
\end{eqnarray}
while the function $Q_{\kappa}$, follows from $R_{\kappa}(x)$, through
\begin{equation}
Q_{\kappa}(x)= \frac{d \,R_{\kappa}(x)}{d\,x} + N_{\kappa}\, x ^{n-1} \ \ .\label{KGGDIII8}
\end{equation}
Eq. (\ref{KGGDIII7}) imposes that $R_{\kappa}$ is a polynomial of degree $n$ i.e.
\begin{equation}
R_{\kappa}(x)=N_{\kappa} \sum_{m=0}^n c_m x^m \ \ , \label{KGGDIII9}
\end{equation}
By direct substitution of the expression of $R_{\kappa}$, in Eq. (\ref{KGGDIII7}), it obtains that the coefficients $c_m = c_m (\kappa)$, of the three highest degree terms of the polynomial $R_{\kappa}$ are given by
\begin{eqnarray}
&&c_n =  \frac{n \kappa^2}{1- n^2 \kappa^2} \ , \label{KGGDIII10} \\
&&c_{n-1}= 0 \ , \label{KGGDIII11} \\
&&c_{n-2}= \frac{n-1}{(1- n^2 \kappa^2)\,[1- (n-2)^2 \kappa^2]} \ ,
 \label{KGGDIII12}
\end{eqnarray}
while the coefficients of order $m$ with $0\leq m \leq n-3$ can be calculated by means of the recursive formula
\begin{eqnarray}
c_m =  \frac{(m+1)(m+2)}{1- m^2 \kappa^2} \,\, c_{m+2}
 \ \ . \label{KGGDIII13}
\end{eqnarray}

After taking into account Eqs. (\ref{KGGDIII8}) and (\ref{KGGDIII9}) it follows that $Q_{\kappa}$ is a polynomial of degree $n-1$, given by
\begin{equation}
Q_{\kappa}(x)= N_{\kappa} \sum_{m=0}^{n-3}(m+1)\, c_{m+1}\, x^{m} + \frac{N_{\kappa}}{1-n^2\kappa^2}\, x^{n-1} \ . \label{KGGDIII14}
\end{equation}

The above described procedure can be used to determine the family of distributions cdf of Type I$_2$ defined through Eq. (\ref{KGGDIII4}). The first member of the family as previously mentioned is just the distribution of Type I$_1$ already discussed. Hereafter are reported the pdf and cdf related to the second and third member of the family:

Second member of the family $(n=2)$
\begin{eqnarray}
p_{\kappa}(x)\!\!\!\!\!\!\!\!\!\!\!&&=  \, (1-4\kappa^2)\,x \exp_{\kappa}(-x),  \\
P_{\kappa}(x)\!\!\!\!\!\!\!\!\!\!\!&&=1\!-\!\left (\! 2\kappa^2x^2\!+\!1 \!+\! x \sqrt{1+\kappa^2x^2}\, \right )\! \exp_{\kappa}(-x) ,  \label{KGGDIII16}
\end{eqnarray}

Third member of the family $(n=3)$
\begin{eqnarray}
p_{\kappa}(x)\!\!\!\!\!\!\!\!\!\!\!&&=\frac{1}{2}\,(1-\kappa^2)\,(1-9\kappa^2)\,x^2 \exp_{\kappa}(-x) , \\
P_{\kappa}(x)\!\!\!\!\!\!\!\!\!\!\!&&=1-\bigg \{\frac{3}{2}\,\kappa^2 (1-\kappa^2)\,x^3 +x \nonumber \\ \!\!\!\!\!\!\!\!\!\!\!&&+ \! \left [ \!1\! + \frac{1}{2}  (1\!-\!\kappa^2)x^2 \right ] \!\! \sqrt{1+\kappa^2x^2} \bigg \}  \exp_{\kappa}(-x) . \label{KGGDIII16}
\end{eqnarray}

\subsection{The distribution of Type I$_9$}

We recall that starting from a given cdf $P_{\kappa}(x)$ with support $x \geq 0$, it is easy to introduce the corresponding cdf $F_{\kappa}(x)$ with support $x\in R$ by means of
$F_{\kappa}(x)=\frac{1}{2}+ \frac{1}{2}\,\frac{x}{|x|}\,P_{\kappa}(|x|)$.

Let us now consider the cdf $F_{\kappa}(x)$ corresponding to the cdf $P_{\kappa}(x)$ of Type I$_9$ i.e. the Half-Normal distribution given by Eq. (\ref{KGGDII5}) with ${\alpha}=2$ and $\nu=1/2)$. Thi corrsponding cdf $F_{\kappa}(x)$ with support $x\in R$ called $\kappa$-deformed Normal cdf, can be written in the form
\begin{equation}
F_{\,\kappa}\,(x)= \frac{1}{2} + \frac{1}{2}\,\,
{\rm erf}_{\kappa}\!\left(\sqrt{\beta}\,x\right) \ \ , \label{KGGDIV2}
\end{equation}
where
\begin{equation}
{\rm
erf}_{\kappa}\!\left(x\right)\!=\!\left(\!1+\frac{1}{2}\,\kappa \!\right)
\!\sqrt{2\kappa}\,\,
\frac{\Gamma\!\left(\frac{1}{2\kappa}+\frac{1}{4}
\right)}{\Gamma\!\left(\frac{1}{2\kappa}-\frac{1}{4}
\right)}\frac{2}{\sqrt{\pi}}\int_0^x \!\exp_{\kappa}(-t^2)\,dt \
\ , \label{KGGDIV3}
\end{equation}
is the $\kappa$-deformed Error function that represents a generalization of the ordinary Error function ${\rm erf}(x)$. The related $\kappa$-deformed Normal pdf $f_{\,\kappa}\,(x)=d F_{\,\kappa}\,(x)/dx$ is an even function with support $x\in R$ and is given by
\begin{equation}
f_{\,\kappa}\,(x)= \sqrt{\frac{2\beta\kappa}
{\pi}}\left(1+\frac{1}{2}\kappa\right)
\frac{\Gamma\!\left(\frac{1}{2\kappa}+\frac{1}{4}
\right)}{\Gamma\!\left(\frac{1}{2\kappa}-\frac{1}{4}
\right)}\,\,\, \exp_{\kappa}(-\beta \,x^2) \ \ . \label{KGGDIV4}
\end{equation}
The moments of $f_{\,\kappa}\,(x)$ of odd order are equal to zero. The variance of $f_{\,\kappa}\,(x)$ has a finite value for $\kappa <2/3$ and assumes the form
\begin{eqnarray}
\sigma^2=\frac{1}{\beta}\,\,\frac{2+\kappa}{2-\kappa}\,\,
\frac{4\kappa}{4-9\kappa^2}\,\Bigg [\frac{
\Gamma\left(\frac{1}{2\kappa}+\frac{1}{4}\right)}
{\Gamma\left(\frac{1}{2\kappa}-\frac{1}{4}\right)}\Bigg ]^2
 \ \ . \ \ \label{KGGDIV5}
\end{eqnarray}

\section{Distribution of Type ${\rm II}$}

Let us consider the cdf $P_{\kappa}(x)$ defined for $x\geq 0$ and given by Eq.(\ref{MSTI0A1}). The positive function $f(x)$ is given by Eq.(\ref{MSTI0A4}) while the positive function $h(x)$ is postulated now to have the following new form
\begin{eqnarray}
h(x)=\frac{f'(x)}{\sqrt{1+\kappa^2\,f(x)^2}} \ \ . \label{MSTI0A5}
\end{eqnarray}

Given the particular structure of the function $h(x)$ which explicit expression becomes
\begin{eqnarray}
h_{\kappa}(x)=\frac{{\alpha} \,\beta\,
x^{{\alpha}-1}}{\sqrt{1+\kappa^2\beta^2x^{2{\alpha}}}} \ , \label{MSTII0A41}
\end{eqnarray}
the two integrals in the formula defining the cdf $P_{\kappa}(x)$, can be performed easily so that $P_{\kappa}(x)$ and the corresponding pdf $p_{\kappa}(x)$, assume the simple forms
\begin{eqnarray}
&&P_{\kappa}(x)=1-\exp_{\kappa}(-\beta
x^{\alpha}) \  ,
 \label{MSTII08}
\\
&&p_{\kappa}(x)=h_{\kappa}(x)\,
\exp_{\kappa}(-\beta\,x^{\alpha}) \ . \label{MSTII09}
\end{eqnarray}
The latter distribution of Type II, represents the $\kappa$-deformed version of the Weibull distribution, which recovers in the $\kappa \rightarrow 0$ limit.

Here after we focus on the expression of the survival or reliability function $S_{\kappa}(x) = 1 - P_{\kappa}(x)$, of the present statistical model in order to better understand the meaning of the function $h_{\kappa}(x)$. It follows that
\begin{eqnarray}
S_\kappa(x)=\exp_\kappa(-\beta\,x^\alpha),  \ .\label{21aa}
\end{eqnarray}
and then
\begin{eqnarray}
p_\kappa(x)= h_\kappa(x)\,S_\kappa(t). \label{28}
\end{eqnarray}
so that $h_{\kappa}(x)$ represents the \emph{hazard function} of the model.

The rate equation for the survival function assumes the form of the first-order linear ordinary differential equation
\begin{eqnarray}
\frac{d S_\kappa(x)}{d x}= -h_\kappa(x)\,S_\kappa(x) \ , \label{27}
\end{eqnarray}
with $S_\kappa(0)=1$. From the latter equation it follows that $S_\kappa(x)$ can be written in the form
\begin{eqnarray}
S_\kappa(x)= e^{-H_\kappa(x)} \ ,\label{33aa}
\end{eqnarray}
where $H_\kappa(x)$ is the cumulative hazard function $H_\kappa(x)$, defined by means of the integral $H_\kappa(x)=\int_0^x h_\kappa(u)\,du$. After performing the integration it obtain the following explicit expression of the cumulative hazard function
\begin{eqnarray}
H_\kappa(x)=\frac{1}{\kappa}\,{\rm arcsinh}\, (\kappa\,\beta\,x^\alpha) \ ,\label{33}
\end{eqnarray}
 which in the $\kappa \rightarrow 0$ limits reduces to the standard Weibull cumulative hazard function $H(x)=\beta\,x^\alpha$.

From Eqs. (\ref{21aa}), (\ref{33aa}) and (\ref{33}) it follows the already known second representation of the $\kappa$-exponential function i.e. $\exp_\kappa(x)=\exp\left(\frac{1}{\kappa} ({\rm arcsinh } (\kappa\,x)\right)$.

The moment of order $m$ of $p_\kappa(x)$ \cite{SupplementaryMaterial}, is given by
\begin{eqnarray}
<x^m>\,=\frac{|2\kappa\beta|^{-m/\alpha}}
{1+\kappa\frac{m}{\alpha}}\,\,\frac{\Gamma\left(\frac{1}{2\kappa}-\frac{m}{2\alpha}\right)}
{\Gamma\left(\frac{1}{2\kappa}+\frac{m}{2\alpha}\right)}\,\,\Gamma\left(1+\frac{m}{\alpha}\right)
 \ \ .  \label{MSTII010}
\end{eqnarray}

Also $G_{\kappa}=1- \int_0^{\infty}[1-P_{\kappa}(x)]^2dx / \int_0^{\infty}[1-P_{\kappa}(x)]\,dx$
i.e. the Gini coefficient can be expressed in closed form
\begin{eqnarray}
G_{\kappa}=1- \frac{\alpha +\kappa}{\alpha+ \frac{1}{2}\kappa}\frac{\Gamma\left(\frac{1}{\kappa}-\frac{1}{2\alpha}\right)}
{\Gamma\left(\frac{1}{\kappa}+\frac{1}{2\alpha}\right)}
\frac{\Gamma\left(\frac{1}{2\kappa}+\frac{1}{2\alpha}\right)}
{\Gamma\left(\frac{1}{2\kappa}-\frac{1}{2\alpha}\right)}
 \ \ .  \label{MSTII010a1}
\end{eqnarray}

The mode of the pdf is located at
\begin{eqnarray}
\!\!\!\!\!\!\!\!\! x_{mode}\,=&&\!\!\!\!\!\!\!\!\!\beta^{-1/\alpha}\left ( \frac{\alpha^2 +2 \kappa^2 (\alpha -1)}{2\kappa^2 (\alpha^2-\kappa^2)}
 \right )^{\!1/2\alpha} \nonumber \\
 &&\!\!\!\!\!\!\!\!\!\times \!\left (\! \sqrt{1+\frac{4\kappa^2 (\alpha^2-\kappa^2)(\alpha -1)^2}{[\alpha^2 + 2 \kappa^2 (\alpha-1)]^2}} -1\!
 \right )^{\!1/2\alpha}
 \!\!\!\! ,  \label{MSTII010a}
\end{eqnarray}
if $\alpha > 1$; for $\alpha=1$ the pdf is a monotonically decreasing function $p_{\kappa}(0)=\beta$; otherwise, the distribution is zero-modal with a pole at the origin.

The quantile function, after inversion of the cdf given by Eq. (\ref{MSTII08}), is available in the following closed form
\begin{eqnarray}
x_{\kappa}(P)=\beta^{-1/\alpha}\left ( \ln_{\kappa} \frac{1}{1 - P} \right )^{1/\alpha}
 \ \ ,  \label{MSTII010b}
\end{eqnarray}
with $0\leq P \leq 1$ and $\ln_{\kappa} t=(t^{\kappa}-t^{-\kappa})/2\kappa$ the $\kappa$-logarithm function. Starting from the expression of the quantile function and after posing $P=1/2$, it obtains the median of the distribution as $x_{med}=\beta^{1/\alpha} \left ( \ln_{\kappa} 2 \right )^{1/\alpha}$.

The asymptotic behavior of $p_{\,\kappa}(x)$, is defined by
\begin{eqnarray}
&&p_{\,\kappa}(x)\,\,
{\atop\stackrel{\textstyle\sim}{\scriptstyle x\rightarrow +\infty}}
\,\,\,\frac{\alpha}{\kappa}\,(2\,\kappa\,\beta)^{-1/\kappa}\,\,
x^{-1-\alpha/\kappa} \ \ , \label{MSTII011} \\
&&p_{\,\kappa}(x)\,\,
{\atop\stackrel{\textstyle\sim}{\scriptstyle x\rightarrow 0^+}}
\,\,\,\alpha\,\beta\,\,\, x^{\alpha-1} \ \ . \label{MSTII012}
\end{eqnarray}

\subsection{Distribution of Type ${\rm II_1}$}

The cdf and the pdf given by Eqs. (\ref{MSTII08}) and (\ref{MSTII09}) in the special case $\alpha=1$ represent a $\kappa$-deformation of the ordinary exponential distribution. In this case the moment of order $m<1/\kappa$ simplifies to
\begin{eqnarray}
<x^m>\,=\frac{\beta^{-m}\,m!}{
\prod_{n=0}^{m}[1-(2n-m)\,\kappa]}\,\,
 \ \ .  \label{MSTII14}
\end{eqnarray}
In particular the mean value and the variance of $p_{\kappa}(x)$ assumes a very simple form and are given by
\begin{eqnarray}
<x>=\frac{1}{\beta}\,\frac{1}{1-\kappa^2} \ \ ; \ \
\sigma^2=\frac{1}{\beta^2}\,
\frac{1+2\kappa^4}{(1-4\kappa^2)(1-\kappa^2)^2} \ \ .
\label{MSTII16}
\end{eqnarray}

The function
${\cal L}_{\kappa}(P)=\int_{\,0}^{x_{\kappa}(P)}\,t\,p(t)
\,dt\,/\int_{\,0}^{\infty}\,t\,p(t)\,dt$ defines the Lorenz curve which assumes the explicite form
\begin{eqnarray}
{\cal L}_{\kappa}(P)=1+\frac{1-\kappa}{2\kappa}(1-P)^{1+\kappa}
 -\frac{1+\kappa}{2\kappa}(1-P)^{1-\kappa} \ \ ,
 \label{MSTII19}
\end{eqnarray}
while the Gini coefficient simplifies to
\begin{eqnarray}
 G_{\kappa}=\frac{2+\kappa^2}{4-\kappa^2} \ \ .
 \label{MSTII19}
\end{eqnarray}

\subsection{Distribution of Type ${\rm II_2}$} Corresponds to $\alpha=2$ and defines the $\kappa$-deformed version of Rayleigh distribution.

\section{Distribution of Type ${\rm III}$}

The cdf of the model of Type III is defined for $x\geq 0$, according to
\begin{equation}
P_{\kappa}(x)\!=\!\frac{\int_{0}^{x}\!h(t)
\exp_{\kappa}[-f(t)]\,dt}{\int_{0}^{x}\!h(t)
\exp_{\kappa}[-f(t)]\,dt \!+\!\lambda \int_{x}^{\infty}\!h(t)
\exp_{\kappa}[-f(t)]\,dt} \ , \label{MSTIII01}
\end{equation}
where $\lambda >0$, while the functions $f(x)$ and $h(x)$ are the same of the model of Type II and are given by Eq.(\ref{MSTI0A4}) and Eq. (\ref{MSTII0A41}) respectively. It results that the model of Type II, emerges to be a special case of the model of Type III, corresponding to $\lambda=1$ called also model of Type III$_1$ .

After performing the integrals in Eq. (\ref{MSTIII01}), the obtaining cdf and the corresponding pdf assume the forms

\begin{eqnarray}
&&P_{\kappa}(x)
=\frac{1- \exp_{\kappa}(-\beta x^{\alpha})}
{1+(\lambda-1)\exp_{\kappa}(-\beta
x^{\alpha})} \label{MSTIII07} \ \ ,
\\
&&p_{\kappa}(x)=\lambda\,h_{\kappa}(x)\,\frac{
\exp_{\kappa}(-\beta x^{\alpha})}
{\big [1+(\lambda-1)\exp_{\kappa}(-\beta
x^{\alpha})\big ]^2} \label{MSTIII08} \ \ , \ \ \ \
\end{eqnarray}
while the survival function $S_{\kappa}(x)\!=\!1\!-\!P_{\kappa}(x)$ is given by
\begin{eqnarray}
S_{\kappa}(x)
=\frac{\lambda}
{\exp_{\kappa}(\beta x^{\alpha})+\lambda-1} \label{MSTIII11C}  \ \ .
\end{eqnarray}

By direct comparison of Eq. (\ref{MSTIII08}) and Eq. (\ref{MSTIII11C})
it follows the relation linking $p_{\kappa}(x)$ and $S_{\kappa}(x)$
\begin{eqnarray}
p_{\kappa}(x)=h_{\kappa}(x)\,S_{\kappa}(x)\left (1- \frac{\lambda -1}{\lambda} S_{\kappa}(x)\right )
 \ \ ,  \label{MSTIII10}
\end{eqnarray}
and finally after recalling that $p_{\kappa}(x)=-d\,S_{\kappa}(x)/dx$ the evolution equation for the $S_{\kappa}(x)$ assume the form
\begin{eqnarray}
\frac{d\,S_{\kappa}(x)}{dx}=- h_{\kappa}(x)\,S_{\kappa}(x)\left (1- \frac{\lambda -1}{\lambda} S_{\kappa}(x)\right )
 \ \ ,  \label{MSTIII11S}
\end{eqnarray}
with $S_{\kappa}(0)=1$. The solution of this equation produces just the expression of  $S_{\kappa}(x)$ as given by Eq. (\ref{MSTIII11C}).

The cumulative hazard function $H_{\kappa}(x)$ is the same of the model of Type II and is given by Eq. (\ref{33}). After performing a variable change the rate equation (\ref{MSTIII11S}) semplifies to $\frac{d S_{\kappa}}{d H_{\kappa}}=-S_{\kappa}\left(1\!-\!\frac{\lambda -1}{\lambda} S_{\kappa} \right )$ with solution given by $S_{\kappa}=\lambda / \left ( e^{H_{\kappa}} + \lambda -1 \right )$. From this rate equation it follows that also the $\kappa$-deformed model of Type III continuous to describe a population kinetics with a bosonic $(0<\lambda <1)$ or fermionic $(\lambda >1)$ character like in the case of the ordinary undeformed model.

\subsection{Distribution of Type ${\rm III_2}$} In the special case where  $\lambda=2$ and $\alpha=1$, the cdf defined for $x\geq 0$ reads
\begin{eqnarray}
P_{\,\kappa}\,(x)=\frac{1-\exp_{\kappa}(-\beta
x)}{1+\exp_{\kappa}(-\beta x)} \ \ , \label{MSTIII11}
\end{eqnarray}
and can be viewed as the $\kappa$-deformed version of the ordinary half-Logistic model which recover in the $\kappa \rightarrow 0$ limit. The survival function as follows from Eq. (\ref{MSTIII11C}), is equal to two times the $\kappa$-deformed Fermi-Dirac function.

The $\kappa$-Logistic cdf defined for $x \in R$ assume the form
\begin{eqnarray}
P_{\,\kappa}\,(x)=\frac{1}{1+\exp_{\kappa}(-\beta
x)} \ \ . \label{MSTIII1S1}
\end{eqnarray}

\section{Distributions of Type ${\rm IV}$}

After recalling the relationship $\exp_{\kappa}(-\beta x^{\alpha})\approx (2\kappa \beta x^{\alpha})^{-1/\kappa}$, with $\alpha>0$, holding for $x \rightarrow \infty$, the following cdf of Type IV
\begin{eqnarray}
P_{\kappa}(x)=(2\kappa\beta)^{1/\kappa}\, \,
x^{\alpha/\kappa}\,\exp_{\kappa}(-\beta
 x^{\alpha}) \ \ ,
 \label{MSTIV1}
\end{eqnarray}
 can be introduced. The related pdf assumes the form
\begin{eqnarray}
p_{\kappa}(x)=\frac{\alpha}{\kappa}\,\,
(2\beta\kappa)^{1/\kappa} \!\left(1-\frac{\kappa \,\beta\,
x^{{\alpha}}}{\sqrt{1+\kappa^2\beta^2x^{2{\alpha}}}}\right) \nonumber \\
\times \,x^{-1+\alpha/\kappa}\, \exp_{\kappa}(-\beta\,x^{\alpha}) \ \ . \label{MSTIV4}
\end{eqnarray}
It is noteworthy that the above distributions do not admit a classic counterpart because in the $\kappa \rightarrow 0$ limit, both $P_{\kappa}(x)$ and $p_{\kappa}(x)$ reduces to zero.

The asymptotic behavior of $P_{\kappa}(x)$ is described by
\begin{eqnarray}
&&P_{\,\kappa}(x)\,\,
{\atop\stackrel{\textstyle\sim}{\scriptstyle x\rightarrow 0^+}}
\,\,\,(2\,\kappa\,\beta)^{1/\kappa}\, x^{\alpha/\kappa} \ \ ,
\label{MSTIIV2}
\\
&&P_{\,\kappa}(x)\,\,
{\atop\stackrel{\textstyle\sim}{\scriptstyle x\rightarrow +\infty}}
\,\,\,1-\frac{ x^{-2\alpha}}{4\kappa^{3}\beta^2}\,\, \ \ .
\label{MSTIV3}
\end{eqnarray}

The moment of order $m<2\alpha$, of $p_{\kappa}(x)$ \cite{SupplementaryMaterial}, is given by
\begin{eqnarray}
<x^m>\,=\frac{(2\kappa\beta)^{-m/\alpha}}{1+\kappa m/2\alpha}\,
\,\frac{\Gamma\left(\frac{1}{\kappa}+\frac{m}{\alpha}\right)\,\Gamma\left(1-\frac{m}{2\alpha}\right)}
{\Gamma\left(\frac{1}{\kappa}+\frac{m}{2\alpha}\right)}
 \ \ .  \label{MSTIV5}
\end{eqnarray}

\section{Distributions of Type V}

Let us consider the derivative $p^{(n)}(x)$, of order $n$,  of
the function $p\,(x)=p^{(0)}(x)$, in following called generatrix function, and suppose that for $0\leq n \leq N$ has the properties i)
$p^{(n)}(0)$ is finite  ii) the function $p^{(n)}(x)/p^{(n)}(0)$ is
positive and monotonically decreasing in $x > 0$, iii)
$\int_0^{+\infty}p^{(n)}(x)\, dx < +\infty$.

We postulate the following distributions of order $n \geq 1$
\begin{eqnarray}
P(n,x)= 1-\frac{p^{(n-1)}(x)}{p^{(n-1)}(0)} \ \ ; \ \ p\,(n,x)= -\frac{p^{(n)}(x)}{p^{(n-1)}(0)} \ \ , \label{MSTIX04}
\end{eqnarray}
where $P(n,x)$ is the cdf and $p\,(n,x)=\frac{d\,P(n,x)}{dx}$ the pdf.

The relationship
\begin{equation}
\int_0^{\infty}\!\!\!x^{m}\,p\,(n,x)\,dx=
-m\frac{p^{(n-2)}(0)}{p^{(n-1)}(0)}\int_0^{\infty}\!\!\!x^{m-1} \,p\,(n-1,x)\,dx \ \ ,
\label{MSTIX05}
\end{equation}
if iterated $m$ times permits to obtain the moment of order $m$ of $p\,(n,x)$, with $0\leq m\leq n-1$, as follows
\begin{eqnarray}
\int_0^{\infty}x^{m}\,p\,(n,x)\,dx=(-1)^m\, m! \,\frac{p^{\,(n-1-m)}(0)}{p^{\,(n-1)}(0)} \ \ .
\label{MSTIX07}
\end{eqnarray}

The procedure described here can be applied to any distribution and
can produce in some case new distributions. In the simplest case of the ordinary exponential distribution this procedure does not produce new distributions.

Here after, by applying this procedure, we introduce the family of distributions of Type V, having as generatrix function the pdf of Type $I_1$. For $x \rightarrow \infty$ the pdf $p_{\,\kappa}(n,x)$ of Type V behave as $p_{\,\kappa}(n,x) \propto x^{-n-\frac{1}{\kappa}}$. In the following we consider the fist three elements of the family. The first element i.e. the distribution of order $n=1$, is just the already discussed distribution of Type $II_1$.

The distributions of order $n=2$ and $n=3$, are new and the related cumulative and probability density functions are given by:
\begin{eqnarray}
P_{\,\kappa}(2,x) = 1- \frac{
\exp_{\kappa}(-\beta\,x)}{\sqrt{1+\kappa^2\beta^2x^{2}}} \ \ , \label{MSTIX21}
\end{eqnarray}
\begin{eqnarray}
p_{\,\kappa}(2,x) =
\left[\frac{1}{1+\kappa^2\beta^2x^2}+\frac{\kappa^2\beta x}
{(1+\kappa^2\beta^2x^2)^{3/2}}\right] \nonumber
\\ \times \,\beta\exp_{\kappa}(-\beta \,x)
\ \ , \label{MSTIX22}
\end{eqnarray}
\begin{eqnarray}
P_{\,\kappa}(3,x) =1-
\left[\frac{1}{1+\kappa^2\beta^2x^2}+\frac{\kappa^2\beta x}
{(1+\kappa^2\beta^2x^2)^{3/2}}\right] \nonumber \\ \times\,\exp_{\kappa}(-\beta \,x) \ \ ,
\label{MSTIX31}
\end{eqnarray}
\begin{eqnarray}
p_{\,\kappa}(3,x) =
 \bigg[\frac{1-\kappa^2}{(1+\kappa^2\beta^2x^2)^{3/2}}
 +\frac{3\kappa^2\beta x}{(1+\kappa^2\beta^2x^2)^{2}}
 \nonumber \\ +\frac{3\kappa^4\beta^2 x^2}{(1+\kappa^2\beta^2x^2)^{5/2}}\bigg]
 \beta\exp_{\kappa}(-\beta \,x) \ \ . \label{MSTIX32}
\end{eqnarray}

The mean value and the variance of the pdf of $p_{\,\kappa}(2,x)$ are respectively $<x>\,= \frac{1}{\beta}$ and $\sigma^2= \frac{1}{\beta^2}\, \frac{1+\kappa^2}{1-\kappa^2}$.
It is remarkable that the mean value and the variance of $p_{\,\kappa}(3,x)$ do not depend on the parameter $\kappa$ and are the same of the ordinary exponential function i.e. $<x>\,=  \frac{1}{\beta}$, $\sigma^2=  \frac{1}{\beta^2}$.

In the $\kappa \rightarrow 0$ limit it obtains $P_0 (n, x)=1- \exp (-\beta x)$ and  $p_{0}(n, x)= \beta \exp (-\beta x)$ holding for  $n\in N$.

\section{Distributions with $\alpha < 0$}

Let us consider the pdf of Type I, indicated here by $p_{\kappa}(\alpha,x)$, with $\alpha >0$. The related inverse pdf $p^{inv}_{\kappa}(\alpha,x)$, is defined through $p^{inv}_{\kappa}( \alpha,x)\,|dx|=p_{\kappa}(\alpha,t) \,|dt|$, where $t=1/x$ and it is easy to verify that $p^{inv}_{\kappa}(\alpha,x)=|p_{\kappa}(-\alpha,x)|$. From this last relationship it follows that the distribution of Type I can be considered also when $\alpha<0$.  Eqs. (\ref{KGGDII5}) and (\ref{KGGDII7}) continue to define the cdf and the pdf respectively for any real not null $\alpha$, as long as the normalization constant $N_{\kappa}$ given by Eq.(\ref{KGGDII6}) is replaced by $|N_{\kappa}|$.

The same analysis can be applied also to the distributions of Type II, III and IV. Alternatively it can be observed that Eqs. (\ref{MSTII08}), (\ref{MSTIII07}) and (\ref{MSTIV1}), defining the three models respectively, in the case $\alpha <0$ define the survival functions $S_{\kappa}(x)$ of the models. The related probability density functions $-\frac{d}{dx}S_{\kappa}(x)$, becomes $|\,p_{\kappa}(x)\,|$ where the functions $p_{\kappa}(x)$ for the three models are given by Eqs. (\ref{MSTII09}), (\ref{MSTIII08}) and (\ref{MSTIV4}) respectively. The parameter $\alpha < 0$ has the same meaning as in the case of ordinary statistics.

\section{Appendix}

The present appendix contains the proofs of the formulas giving the Mellin transform of $\exp_{\kappa}(-x)$, the normalization constant of the pdf of Type I, the moment of order $m$ of the pdf of Type I, the moment of order $m$ of the pdf of Type II and the moment of order $m$ of the pdf of Type IV. The present appendix appears as "Supplementary Material" in the website of the journal Europhysics Letters (EPL) together with the present paper.

\subsection{Mellin transform of $\exp_{\kappa}(-x)$}

 Here it is obtained the Mellin transform of the function $\exp_{\kappa}(-x)=(\sqrt{1+\kappa^2 x^2}-\kappa x)^{1/\kappa}$ i.e. its moment of order $r-1$
\begin{eqnarray}
{\cal M}_{\kappa}(r)=\int_0^{\infty} x^{r-1}\exp_{\kappa}(-x)\, dx
 \ \ .  \label{I1}
\end{eqnarray}

After performing a variable change through the transformation $w=\left( \sqrt{1+\kappa^2 x^2}-\kappa x\right)^2$ from which it follows $x=\frac{1}{2\kappa}w^{-\frac{1}{2}}(1-w)$, $dx=-\frac{1}{4\kappa}w^{-\frac{3}{2}}(1+w)\, dw$ and $\exp_{\kappa}(-x)=w^{\frac{1}{2\kappa}}$.

The integration can be performed as follows
\begin{eqnarray}
{\cal M}_{\kappa}(r)\!\!\!\!\!\!\!\!\!&&=\frac{1}{4\kappa}\, (2 \kappa)^{1-r} \!
\int_0^{1} {w}^{ \frac{1}{2\kappa}-\frac{r}{2}-1} \left (1-w \right )^{r -1} \left (1+w \right ) \, dw \nonumber \\
&&=\frac{1}{2}\, (2 \kappa)^{-r} \! \int_0^{1} {w}^{ \frac{1}{2\kappa}-\frac{r}{2}-1} \left (1-w \right )^{r -1} \, dw \nonumber \\
&&+ \, \frac{1}{2}\, (2 \kappa)^{-r} \! \int_0^{1} {w}^{ \frac{1}{2\kappa}-\frac{r}{2}} \left (1-w \right )^{r -1} \, dw \nonumber \\
&&=\frac{1}{2}\, (2 \kappa)^{-r} \, B \left( \frac{1}{2\kappa}-\frac{r}{2}\, \,  , \, \, r  \right )
\nonumber \\
&&+ \, \frac{1}{2}\, (2 \kappa)^{-r} \, B \left( \frac{1}{2\kappa}-\frac{r}{2}+ 1 \, \,  , \, \, r  \right ) \nonumber \\
&&=\frac{1}{2}\, (2 \kappa)^{-r} \, \,
\frac{\Gamma \left ( \frac{1}{2\kappa}-\frac{r}{2}\right)\, \Gamma \left ( r \right) }{\Gamma \left (  \frac{1}{2\kappa}+\frac{r}{2} \right)}
\nonumber \\
&&+\, \frac{1}{2}\, (2 \kappa)^{-r} \, \,
\frac{\Gamma \left ( \frac{1}{2\kappa}-\frac{r}{2} +1 \right)\, \Gamma \left ( r \right) }{\Gamma \left (  \frac{1}{2\kappa}+\frac{r}{2} + 1 \right)} \nonumber \\
&&=\frac{1}{2}\, (2 \kappa)^{-r} \, \left ( 1+ \frac{\frac{1}{2\kappa}-\frac{r}{2}}{\frac{1}{2\kappa}+\frac{r}{2}} \right )
\frac{\Gamma \left ( \frac{1}{2\kappa}-\frac{r}{2}\right)
 }{\Gamma \left (  \frac{1}{2\kappa}+\frac{r}{2} \right)} \,\Gamma \left ( r \right)
  ,  \label{I2}
\end{eqnarray}
obtaining in this way the Mellin transform of $\exp_{\kappa}(-x)$ in the form
\begin{eqnarray}
{\cal M}_{\kappa}(r)= \frac{(2 \kappa)^{-r}}{1+\kappa r} \,\,
\frac{\Gamma \left ( \frac{1}{2\kappa}-\frac{r}{2}\right)
 }{\Gamma \left (  \frac{1}{2\kappa}+\frac{r}{2} \right)} \,\, \Gamma \left ( r \right)
 \ \ .  \label{I3}
\end{eqnarray}
holding for $r<1/\kappa$. In the $\kappa \rightarrow 0$ limit, after taking into account
the relationship $\Gamma\left(z-c\right)/\,\Gamma\left(z+c\right) \approx z^{-2c}$
holding in the Stirling approximation for $z\rightarrow \infty$, it obtains ${\cal M}_{0}(r)= \Gamma (r)$.

\subsection{Normalization constant of the the distribution of Type I}

The normalization constant $N_{\kappa}$ appearing in the expression of the pdf
\begin{eqnarray}
p_{\,\kappa}\,(x)= N_{\kappa}
\,\, x^{\alpha\nu-1} \,\exp_{\kappa}(-\beta x^{\alpha}) \ \ , \label{I4}
\end{eqnarray}
follows from the condition $\int_0^{\infty}p_{\kappa}(x)=1$ according to
\begin{eqnarray}
\frac{1}{N_{\kappa}}=\int_0^{\infty} x^{\alpha\nu-1} \,\exp_{\kappa}(-\beta x^{\alpha}) \, dx
 \ \ .  \label{I5}
\end{eqnarray}

After introducing the variable change through $w=\beta x^{\alpha}$ from which follows that $x=\beta^{-\frac{1}{\alpha}} w^{\frac{1}{\alpha}}$ and $dx=\alpha^{-1} \beta^{-\frac{1}{\alpha}} w^{\frac{1}{\alpha}-1} \, dw$, the integral defining the normalization constant becomes
\begin{eqnarray}
\frac{1}{N_{\kappa}}&&\!\!\!\!\!= \frac{1}{\alpha \beta^{\,\nu}} \int_0^{\infty} w^{\nu-1} \,\exp_{\kappa}(-w) \, dw
\noindent \\
&&\!\!\!\!\!= \frac{1}{\alpha \beta^{\,\nu}} {\cal M}_{\kappa}(\nu)
 \ \ .  \label{I5}
\end{eqnarray}
After takining into account the expression of ${\cal M}_{\kappa}(\nu)$ it follows
\begin{eqnarray}
N_{\kappa}=(1+\kappa \,\nu)\,(2\,\kappa)^{\nu}
 \,\,
\frac{\Gamma\!\left(\frac{1}{2\kappa}+\frac{\nu}{2}
\right)}{\Gamma\!\left(\frac{1}{2\kappa}-\frac{\nu}{2}
\right)}\,\, \frac{\alpha\,\beta^{\,\nu}}{\Gamma\!\left(\nu\right)} \ \ ,\label{KGGDII6}
\end{eqnarray}
with $\nu<1/\kappa$.

\subsection{Moment of the distribution of Type I}

The moment of order $m$,  $<x^m>\,=\int_0^{\infty} x^m \, p_{\kappa}(x)\, dx$
for the pdf of Type I
\begin{eqnarray}
p_{\,\kappa}\,(x)= N_{\kappa}
\,\, x^{\alpha\nu-1} \,\exp_{\kappa}(-\beta x^{\alpha})\, dx \ \ , \label{KGGDII7}
\end{eqnarray}
is defined according to
\begin{eqnarray}
<x^m>\,=N_{\kappa}
\int_0^{\infty} x^{m+\alpha\nu-1} \,\exp_{\kappa}(-\beta x^{\alpha})
 \ \ .  \label{IV3}
\end{eqnarray}
After introducing the variable change $w=\beta x^{\alpha}$, the moment assume form
\begin{eqnarray}
<x^m>\,\,=&&\!\!\!\!\!N_{\kappa} \,\, \frac{1}{\alpha \beta^{\,\nu}} \int_0^{\infty} w^{\frac{m}{\alpha} + \nu-1} \,\exp_{\kappa}(-w) \, dw \nonumber \\
=&&\!\!\!\!\!N_{\kappa} \,\, \frac{1}{\alpha \beta^{\,\nu}} \,\, {\cal M}_{\kappa}\left (\nu +\frac{m}{\alpha} \right) \ \ .  \label{I5}
\end{eqnarray}

Finally by taking into account the expressions of the normalization constant $N_{\kappa}$ and the Mellin transform ${\cal M}_{\kappa}\left (\nu +\frac{m}{\alpha} \right)$, obtained previously, the moment of order $m$ of $p_{\,\kappa}\,(x)$ assumes the form
\begin{eqnarray}
<x^m>&&\!\!\!\!\!\!\!\!=(2\kappa \beta)^{-\frac{m}{\alpha}}\,\frac{1+\kappa\nu}
{1+\kappa\left(\nu+\frac{m}{\alpha}\right)}\,\frac{\Gamma\left(\frac{1}{2\kappa}+\frac{\nu}{2}\right)}
{\Gamma\left(\frac{1}{2\kappa}-\frac{\nu}{2}\right)} \nonumber \\
&&\times
\frac{\Gamma\left(\frac{1}{2\kappa}-\frac{\nu}{2}-\frac{m}{2\alpha}\right)}
{\Gamma\left(\frac{1}{2\kappa}+\frac{\nu}{2}+\frac{m}{2\alpha}\right)}\,\frac{\Gamma\left(\nu+\frac{m}{\alpha}\right)}{\Gamma\left(\nu\right)}
\ \ ,  \label{KGGDII11M}
\end{eqnarray}
with $\nu +m/\alpha<1/\kappa$.

\subsection{Moment of the distribution of Type II}

Let us observe that the pdf of Type II, can be written as
\begin{eqnarray}
p_{\,\kappa}\,(x)= - \frac{d \, S_{\,\kappa}\,(x)}{d\, x}\,  \ \ , \label{II1}
\end{eqnarray}
where the survival function $S_{\kappa}(x)\, dx$ has the expression
\begin{eqnarray}
S_\kappa(x)=\exp_\kappa(-\beta\,x^\alpha),  \ .\label{II2}
\end{eqnarray}

The moment of order $m$ of  $p_{\kappa}(x)\, dx$ can be written as
\begin{eqnarray}
<x^m>&&\!\!\!\!\!=\int_0^{\infty} x^m \, p_{\kappa}(x)\, dx  \nonumber \\
&&\!\!\!\!\!=-\int_0^{\infty} x^m \, \frac{d\,S_{\kappa}(x)}{dx}\, dx  \nonumber \\
&&\!\!\!\!\!=m\int_0^{\infty} x^{m-1} \, S_{\kappa}(x)\, dx  \ \ ,\label{II2}
\end{eqnarray}
where in the last step an integration by part is performed.

After substitution of the expression of $S_{\kappa}(x)\, dx$ in the last integral and after making the variable change $w=\beta x^{\alpha}$ the moment of order $m$ of the pdf $p_{\kappa}(x)\, dx$ becomes
\begin{eqnarray}
<x^m>\,=&&\!\!\!\!\!\!\!\! \frac{m}{\alpha \beta^{\, m/\alpha}} \int_0^{\infty} w^{\frac{m}{\alpha}-1} \,\exp_{\kappa}(-w) \, dw
\noindent \\
=&&\!\!\!\!\!\!\!\! \frac{m}{\alpha \beta^{\, m/\alpha}} \,\, {\cal M}_{\kappa}\left(\frac{m}{\alpha} \right)
 \ \ .  \label{I5}
\end{eqnarray}
Finally by taking into account the expression of Mellin transform of $\exp_{\kappa}(-x)$, the moment of $p_{\kappa}(x)$ assumes the form
\begin{eqnarray}
<x^m>\,=\frac{|2\kappa\beta|^{-\frac{m}{\alpha}}}
{1+\kappa\frac{m}{\alpha}}\,\,\frac{\Gamma\left(\frac{1}{2\kappa}-\frac{m}{2\alpha}\right)}
{\Gamma\left(\frac{1}{2\kappa}+\frac{m}{2\alpha}\right)}\,\,\Gamma\left(1+\frac{m}{\alpha}\right)
 \ \ ,  \label{MSTII010}
\end{eqnarray}
with $m<\alpha/\kappa$.

\subsection{Moment of the distribution of Type IV}

Let us consider the moment of order $m$, defined through the relationship $<x^m>\,=\int_0^{\infty} x^m \, p_{\kappa}(x)\, dx$, of the pdf
\begin{eqnarray}
p_{\kappa}(x)=\frac{\alpha}{\kappa}\,\,
(2\beta\kappa)^{\frac{1}{\kappa}} \!\left(1-\frac{\kappa \,\beta\,
x^{{\alpha}}}{\sqrt{1+\kappa^2\beta^2x^{2{\alpha}}}}\right)
x^{\frac{\alpha}{\kappa}-1} \nonumber \\ \times \exp_{\kappa}(-\beta\,x^{\alpha}) \ \ . \label{IV2}
\end{eqnarray}

After substitution of the expression of $p_{\kappa}(x)$ in the integral defining the moment and after performing the variable change $t=\beta x^{\alpha}$ it obtains
\begin{eqnarray}
<x^m>\,=2^{\frac{1}{\kappa}}\, {\kappa}^{\frac{1}{\kappa}-1}\, {\beta}^{-\frac{m}{\alpha}}
&&\!\!\!\!\!\!\!\!\!\int_0^{\infty} {t}^{\frac{m}{\alpha} + \frac{1}{\kappa}-1}
\left(1-\frac{\kappa \,t}{\sqrt{1+\kappa^2 t^2}} \right)  \nonumber \\ &&\times \exp_{\kappa}(-t) \, dt
 \ \ .  \label{IV3}
\end{eqnarray}

A farther variable change is introduced at this point through $w=\left( \sqrt{1+\kappa^2 t^2}-\kappa t\right)^2$. From this transformation it obtains $t=\frac{1}{2\kappa}w^{-\frac{1}{2}}(1-w)$, $dt=-\frac{1}{4\kappa}w^{-\frac{3}{2}}(1+w)\, dw$, $\sqrt{1+\kappa^2 t^2}=\frac{1}{2} w^{-\frac{1}{2}}(1+w)$ and $\exp_{\kappa}(-t)=w^{\frac{1}{2\kappa}}$.

Therefore the moment assumes the form
\begin{eqnarray}
<x^m>\!\!\!\!\!\!\!\!\!\!&&=\frac{1}{\kappa}\, (2 \kappa \beta )^{-\frac{m}{\alpha}}
\!\int_0^{1} \!{w}^{\left(1-\frac{m}{2\alpha}\right) -1} \left (1-w \right )^{\left(\frac{1}{\kappa}+\frac{m}{\alpha}\right) -1} \, dw \nonumber \\
&&=\frac{1}{\kappa}\, (2 \kappa \beta )^{-\frac{m}{\alpha}} \, B \left( 1-\frac{m}{2\alpha}\, , \,\, \frac{1}{\kappa}+\frac{m}{\alpha}  \right ) \nonumber \\
&&=\frac{1}{\kappa}\, (2 \kappa \beta )^{-\frac{m}{\alpha}} \,
\frac{\Gamma \left ( \frac{1}{\kappa}+\frac{m}{\alpha} \right)\, \Gamma \left ( 1 - \frac{m}{2\alpha} \right) }{\Gamma \left ( 1+ \frac{1}{\kappa}+\frac{m}{2\alpha} \right)}
 \ \ .  \label{IV4}
\end{eqnarray}

Taking into account the property $\Gamma (1+x)= x \,\Gamma (\mathrm{}x)$,  the moment of $p_{\kappa}(x)$ can be written finally in the form
\begin{eqnarray}
<x^m>\,=\frac{(2\kappa\beta)^{-\frac{m}{\alpha}}}{1+\frac{\kappa m}{2\alpha}}\,
\,\frac{\Gamma\left(\frac{1}{\kappa}+\frac{m}{\alpha}\right)\,\Gamma\left(1-\frac{m}{2\alpha}\right)}
{\Gamma\left(\frac{1}{\kappa}+\frac{m}{2\alpha}\right)}
 \ \ ,  \label{IV5}
\end{eqnarray}
with $m<2\alpha$.

\end{document}